\documentclass[11pt]{amsart}
\usepackage{geometry}                
\usepackage{graphicx}
\usepackage{amssymb}
\usepackage[usenames, dvipsnames]{xcolor}
\newcommand{\E}{\mathbb{E}}
\renewcommand{\P}{\mathbb{P}}

\newcommand{\BBr}{\mathbb{R}}
\newcommand{\Var}{\operatorname{Var}}
\newcommand{\Cov}{\operatorname{Cov}}

\newcommand{\Xset}{\mathbb{X}}

\newcommand{\x}{\mathbf{x}}
\newcommand{\X}{\mathbf{X}}

\newcommand{\A}{\mathcal{A}}
\newcommand{\Nlaw}{\mathcal{N}}
\newcommand{\Xmc}{\mathbf{X}_\text{MC}}
\usepackage{xcolor}

\newtheorem{lemma}{Lemma}[section]
\newtheorem{proposition}{Proposition}[section]
\newtheorem{remark}{Remark}[section]

\title{Sequential design of experiments for estimating percentiles of black-box functions}
\author{T. Labopin-Richard \and V. Picheny}
\address{TLR in with the Institut de Math\'ematiques de Toulouse (CNRS UMR 5219). Universit\'e Paul Sabatier, 118 route de Narbonne, 31062 Toulouse, France. VP is with MIAT, Universit\'e de Toulouse, INRA, Castanet-Tolosan, France}
\date{\today}                                 
\begin{document}
\begin{abstract}
Estimating percentiles of black-box deterministic functions with random inputs is a challenging task 
when the number of function evaluations is severely restricted, which is typical for computer experiments.
This article proposes two new sequential Bayesian methods for percentile estimation based on the Gaussian
Process metamodel. Both rely on the Stepwise Uncertainty Reduction paradigm, hence aim at 
providing a sequence of function evaluations that reduces an uncertainty measure associated with
the percentile estimator. The proposed strategies are tested on several numerical examples,
showing that accurate estimators can be obtained using only a small number of functions evaluations.
\end{abstract}
\maketitle
\section{Introduction}\label{sec:introduction}

In the last decades, the question of designing experiments for the efficient exploration and analysis of numerical black-box models 
has received a wide interest, and metamodel-based strategies have been shown to offer efficient alternatives in many contexts,
such as optimization or uncertainty quantification.
We consider here the question of estimating percentiles of the output of a black-box model, with the help of Gaussian Process (GP) metamodels
and sequential sampling. More precisely, let $g: \mathbb{X} \subset \BBr^d \rightarrow \BBr$ denote the output of interest of the model, 
the inputs of which can vary within $\mathbb{X}$.
We assume here that the multivariate input $X$ is modelled as a random vector; then, our
objective is to estimate a percentile of $g(X)$: 
\begin{equation}\label{eq:qdef}
 q_{\alpha}(g(X))=q_{\alpha}(Y)=F_Y^{-1}(\alpha),
\end{equation}
for a fixed level $\alpha \in ]0, 1[$, where $F_U^{-1}:= \inf \{\x : F_U(\x) \geq u \}$ 
denotes the generalized inverse of the cumulative distribution function of a random variable $U$.
We consider here only random vectors $X$ and functions $g$ regular enough to have $F_Y \left( F_Y^{-1}(\alpha)\right)=\alpha$ (that is, $F_Y$ is continuous). Since the level $\alpha$ is fixed, we omit the index in the sequel.

A natural idea to estimate a percentile consists in using its empirical estimator:
having at hand a sample $(X_i)_{i=1, \dots, n}$ of the input law $X$, we run it through the computer model
to obtain a sample $(Y_i)_{i=1, \dots, n}$ of the output $Y$. Then, denoted $Y_{(k)}$ the $k$-th order statistic of the previous sample, the estimator
\begin{equation}\label{eq:qemp}
  q_n:= Y_{(\lfloor n \alpha \rfloor +1)},
\end{equation}
is consistent and asymptotically Gaussian under weak assumptions on the model (see \cite{order} for more details). However, for computationally expensive models, 
the sample size is drastically limited, which makes the estimator (\ref{eq:qemp}) impractical.
In that case, one may replace the sample $(X_i)_i$ by a sequence of well-chosen points that
provide a useful information for the percentile estimation. Besides, if the points $X_i$ are not sampled with the distribution of $X$, 
the empirical percentile (\ref{eq:qemp}) is biased, so another estimator must be used.

In \cite{cannamela}, the authors proposed an estimator based on the variance reduction or on the controlled stratification and give asymptotic results. 
Nevertheless, the most usual approach of this problem is a Bayesian method which consists in assuming that $g$ is the realization of a well-chosen Gaussian process. 
In this context, \cite{oakley} propose a two-step strategy: first, generate an initial set of observations to train a 
GP model and obtain a first estimator of the percentile, then increase the set of observations by a second set likely to improve the estimator. 
In \cite{jala2}, the authors proposed a sequential method (called GPQE and GPQE+ algorithms), based on the GP-UCB algorithm of \cite{regret},
that is, making use of the confidence bounds provided by the Gaussian Process model.

In this paper we propose two new algorithms 
based on \textit{Stepwise Uncertainty Reduction} (SUR),
a framework that has been successfully applied to closely related problem such as optimization \cite{victorMO}, or the dual problem of 
percentile estimation, the estimation of a probability of exceedance \cite{bect,chevalierExcursion}.
A first strategy has been proposed for the percentile case in \cite{flooding} and \cite{jala1} that rely on expensive simulation procedures. 
Nevertheless, finding a statistically sound algorithm with a reasonable cost of computation, in particular when the problem dimension increases, is still an open problem. 

The rest of the paper is organized as follow. In Section 2, 
we introduce the basics of Gaussian Process modelling, 
our percentile estimator and the Stepwise Uncertainty Reduction framework. 
Section 3 
describes our two algorithms to estimate a percentile. 
Some numerical simulations to test the two methods are presented in Section 4, 
followed by concluding comments in Section 5. 
Most of the proofs are deferred to the Appendix.
\section{Gaussian Process model}\label{sec:GP}

\noindent {\bf 2.1. Model definition}\label{sec:GPmodel}

We consider here the classical Gaussian Process framework in computer experiments \cite{sacks1989design,covfunction,rasmussen2006gaussian}: we suppose that $g$ is the realization of a GP denoted by $G(.)$ with known
mean $\mu$ and covariance function $c$. 

Given an observed sample  $\A_n=\{(\x_1, g_1), (\x_2, g_2), \dots (\x_n,g_n)\}$  with all $\x_i \in \Xset$ and $g_i=g(x_i)$, the distribution of $G | \A_n$ is entirely known:
\begin{equation*}
\mathcal{L}\left(G| \mathcal{A}_n)= GP(m_n(.), k_n(., .)\right),
\end{equation*}
with, $\forall \x \in \mathbb{X}$,
\begin{equation*}
\begin{aligned}
m_n(\x)&=\E(G(\x) | \mathcal{A}_n)= c_n(\x)^TC_n^{-1}\mathbf{g_n},\\
k_n(\x, \x')&=\Cov\left(G(\x), G(\x') | \mathcal{A}_n \right)=c(\x, \x')-c_n(\x)^TC_n^{-1}c_n(\x'),\\
\end{aligned}
\end{equation*}
where we denote $c_n(\x)=[c(\x_1, \x), \dots, c(\x_n, \x)]^T$, $C_n=[c(\x_i, \x_j)]_{1 \leq i, j \leq n}$ and $\mathbf{g_n}=[g_1, \dots, g_n]$. In the sequel, we also denote $s_n^2(\x)=k_n(\x, \x)$. 

We use here the standard \textit{Universal Kriging} framework, 
where the covariance function depends on unknown parameters that are inferred from $\A_n$, 
using maximum likelihood estimates for instance. Usually, the estimates are used as face value,
but updated when new observations are added to the model.

In the sequel, we will need the following property of the kriging model (see \cite{MAJkrigeage}).
\begin{proposition}
\label{MAJkrigeage}
Moments at step $n+1$ are linked to the moments at step $n$ by the one-step update formula:
\begin{equation*}
\left\{
\begin{aligned}
m_{n+1}(\x)&= m_n(\x) + \frac{k_n(\x_{n+1}, \x)}{s_n^2(\x_{n+1})}\left(g_{n+1} - m_n(\x_{n+1}) \right)\\
s^2_{n+1}(\x)&= s_n^2(\x) - \frac{k_n^2(\x_{n+1}, \x)}{s_n^2(\x_{n+1})}\\
k_{n+1}(\x, \x')&= k_n(\x, \x') - \frac{k_n(\x_{n+1},\x) k_n(\x_{n+1}, \x')}{s_n(\x_{n+1})^2}\\
\end{aligned}
\right.,
\end{equation*}
where $(\x_{n+1}, g_{n+1})$ is a new observational event.
\end{proposition}
%

\noindent {\bf 2.2 Percentile estimation}\label{sec:GPpercentile}

Since each call to the code $g$ is expensive, the sequence of inputs to evaluate, $\{\x_1, \ldots, \x_n\}$, must be chosen carefully
to make our estimator as accurate as possible.
The general scheme based on GP modelling is of the following form:
\begin{itemize}
\item[$\bullet$] For an initial budget $N_0$, we build an initialisation sample $(\x_0^i, g(\x_0^i))_{i=1...N_0}$, typically using a space-filling strategy, and compute the estimator of the percentile $q_{N_0}$. 
\item[$\bullet$] At each step $n+1 \geq N_0+1$ and until the budget $N$ of evaluations is reached: knowing the current set of observations $\mathcal{A}_{n}$ and estimator $q_{n}$, we choose the next point to evaluate $\x_{n+1}^*$, based on a so-called \textit{infill criterion}. We evaluate $g(\x_{n+1}^*)$ and update the observations $\mathcal{A}_{n+1}$ and the estimator $q_{n+1}$. 
\item[$\bullet$] $q_N$ is the estimator of the percentile to return.  
\end{itemize}
In the following, we describe first the estimator we choose, then the sequential strategy adapted to the estimator.

\noindent {\bf Percentile estimator.}
First, from a GP model we extract a percentile estimator. Considering that, conditionally on $\A_n$, the best approximation of $G(\x)$ is $m_n(\x)$,
an intuitive estimator is simply the percentile of the GP mean:
\begin{equation}\label{eq:estimateur1}
q_n:=q(m_n(X)) = q \left( \E \left[G(X) | \mathcal{A}_n \right] \right).
\end{equation}
This is the estimator chosen for instance in \cite{oakley}.

Another natural idea can be to consider the estimator
that minimizes the mean square error $\E \left( (q - {q}_n )^2 \right)$ among all $\mathcal{A}_n$-measurable estimator:
\begin{equation}\label{eq:estimateur2}
q_{n} = \E \left( q({G(X)}) | \mathcal{A}_n \right).
\end{equation}
This estimator is used for instance in \cite{jala2}. 
Despite its theoretical qualities, this estimator suffers from a major drawback, 
as it cannot be expressed in a computationally tractable form, and must be estimated using simulation techniques,
by drawing several trajectories of $G(.)$, computing the percentile of each trajectory and averaging.

Hence, in the sequel, we focus on the estimator (\ref{eq:estimateur1}), which allows us to derive
closed-form expressions of its update when new observations are obtained,
as we show in the next section. 

\begin{remark}
In the case of the dual problem of the probability of failure estimation $u(g)=\P(g(X) >u)$, this later estimator is easier to compute. 
Indeed, is shown in \cite{bect}:
\begin{equation*}
 \E\left( u(G) | \mathcal{A}_n \right) = \E \left( \int_{\mathbb{X}} \mathbf{1}_{G>u} dP_X \right) = \int_{\mathbb{X}} p_n dP_X,
\end{equation*}
where $p_n(\x)= \Phi \left( \frac{m_n(\x)-u}{s_n(\x)} \right)$, for $\Phi $ the cumulative distribution function of the standard Gaussian distribution. 
This compact form is due to the possibility to swap the integral and expectation, which is not feasible for the percentile case. 
\end{remark} 

\noindent {\bf Sequential sampling and Stepwise Uncertainty Reduction.}
We focus here on methods based on the sequential maximization of an infill criterion, that is, of the form:
\begin{equation}
\label{critere}
\x_{n+1}^*=\underset{\x_{n+1} \in \mathbb{X}}{\operatorname{argmax}\, } J_{n}(\x_{n+1}),
\end{equation}
where $J_{n}$ is a function that depends on $\mathcal{A}_{n}$ (through the GP conditional distribution) and $q_{n}$. 

Intuitively, an efficient strategy would explore $\mathbb{X}$ enough to obtain a GP model reasonably accurate everywhere, 
but also exploit previous results to identify the area with response values close to the percentile and sample more densely there.

To this end, the concept of Stepwise Uncertainty Reduction (SUR) has been proposed originally in \cite{geman1996active} 
as a trade-off between exploitation and exploration, and has been successfully adapted to optimization \cite{villemonteix2009informational,victorMO}
or probability of failure estimation frameworks \cite{bect,chevalierExcursion}.
The general principle of SUR strategies is to define an uncertainty measure related to the objective pursued,
and add sequentially the observation that will reduce the most this uncertainty.
The main difficulty of such an approach is to evaluate the potential impact of a candidate point $\x_{n+1}$ without having access to $g(\x_{n+1})=g_{n+1}$ 
(that would require running the computer code).

In the percentile estimation context, \cite{jala1} and \cite{flooding} proposed to choose the next point to evaluate as the minimizer of the conditional variance of the percentile estimator (\ref{eq:estimateur2}).
%
This strategy showed promising results, as it substantially outperforms pure exploration, and, in small dimension, 
manages to identify the percentile area (that is, where $f$ is close to its percentile) and choose the majority of the points in it. 
Nevertheless, computing $V_{n}(\x)$ in \cite{jala1} or \cite{flooding} is very costly, as it requires drawing many GP realizations, which hinders its use in practice
for dimensions larger than two.

In next section, we propose other functions $J$ that also quantify the uncertainty associated with our estimator, but which have closed forms and then are less expensive to compute. 
To do so, we first exhibit the formula to update the current estimator $q_{n}$ (build on $\mathcal{A}_{n}$) to $q_{n+1}(\x_{n+1})$ the estimator at step $n+1$ if we had chosen $\x_{n+1}^*=\x_{n+1}$.

\section{Main results}\label{sec:algo}

\noindent {\bf 3.1. Update formula for the percentile estimator}\label{sec:updateestim}

In this section, we express the estimator $q_{n+1}(\x_{n+1})$ as a function of the past observations $\mathcal{A}_n$, the past percentile estimator $q_n$, a candidate point $\x_{n+1}$ and its evaluation $g_{n+1}$. 

We focus on the estimator ($\ref{eq:estimateur1})$, which is at step $n$ the percentile of the random vector $m_n(X)$. 
Since no closed-form expression is available, we approach it by using the empirical percentile. 
Let $\Xmc=(\x_{\text{MC}}^1, \dots, \x_{\text{MC}}^l)$ be an independent sample of size $l$, distributed as $X$. We compute $m_n(\Xmc)$ and order this vector by denoting $m_n(\Xmc)_{(i)}$ the $i$-th coordinate. Then we choose
\begin{equation}\label{eq:estimator}
q_n=m_n(\Xmc)_{(\lfloor l \alpha \rfloor +1)}.
\end{equation} 

\begin{remark}
Since the observation points $(\x_1, \dots,\x_n)$ do not follow the distribution of $X$, they cannot be used to estimate the percentile.
Hence, a different set $(\Xmc)$ must be used.
\end{remark}

In the sequel, we denote by $\x_n^q$ the point of $\Xmc$ such that
\begin{equation*}
\begin{aligned}
q_n=m_n(\x_n^q).
\end{aligned}
\end{equation*}
This point is referred to as \textit{percentile point}.

Now, let us consider that a new observation $g_{n+1} = g(\x_{n+1})$ is added to $\A_n$.
The key of SUR strategies is to measure the impact of this observation on our estimator $q_n$, that is, express 
$q_{n+1}=m_{n+1}(\x_{n+1}^q)$ as a function of $g_{n+1}$ and $\x_{n+1}$. 

First, by Proposition \ref{MAJkrigeage}, we have:
\begin{equation}
m_{n+1}(\Xmc)=m_n(\Xmc) + \frac{k_n(\Xmc, \x_{n+1})}{s_n(\x_{n+1})^2}\left(g_{n+1}-m_n(\x_{n+1}) \right).
\end{equation}
We see directly that once $\x_{n+1}$ is fixed, all the vector $m_{n+1}(\Xmc)$ is determined by the value of $g_{n+1}$.
Our problem is then to find, for any $g_{n+1}$ in $\mathbb{R}$, which point of $\Xmc$ is the percentile point, that is, 
which point satisfies
\begin{equation}
 m_{n+1}(\Xmc)_{\lfloor l \alpha \rfloor +1}=m_{n+1}\left(\x_{n+1}^q \right).
\end{equation}

Let us denote 
$\mathbf{b}=m_n(\Xmc)$ and $\mathbf{a}=k_n(\Xmc, \x_{n+1})$, which are vectors of $\mathbb{R}^l$, 
and $z=\frac{g_{n+1}-m_n(\x_{n+1})}{s_n^2(\x_{n+1})}$, so that the updated mean simply writes 
as a linear function of $z$, $\mathbf{b} + \mathbf{a}z$. Our problem can then be interpreted 
graphically: each coordinate of $m_{n+1}(\Xmc)$ is represented by a straight line of equation:
\begin{equation}
 b_i+a_iz, \, \, i \in \{1, \dots l\},
\end{equation}
and the task of finding $\x_{n+1}^q$ for any value of $g_{n+1}$ amounts to finding the 
$\lfloor l \alpha \rfloor +1$ lowest line for any value of $z$.


We can first notice that the lines order changes only when two lines intersect each other. 
There are at most $L=\frac{l(l-1)}{2}$ intersection points, which we denote $I_1, \ldots I_L$, in increasing order. 
We set $I_0=- \infty$ and $I_{L+1}= + \infty$, and introduce $(B_i)_{0 \leq i \leq L}$, the sequence of intervals between intersection points:
\begin{equation} \label{Bi}
B_i = [I_i, I_{i+1}]  \, \,  \text{for} \, \, i \in \left[0, L\right]\\
\end{equation}
For any $z \in B_i$, the order of $(b_i+a_i z)$ is fixed. 

Denoting $j_i$ the index of the $\lfloor l \alpha \rfloor  + 1 $ lowest line, 
we have:
\begin{equation}\label{pointPercentile}
  \x_{n+1}^q = \x_{\text{MC}}^{j_i}, \quad z \in B_i,
\end{equation}
the percentile point when $z \in B_i$, which we henceforth denote $\x_{n+1}^q(B_i)$.

Finally, we have shown that

\begin{proposition}\label{MAJpercentile}
Under previous notations, at step $n$ (when we know $\mathcal{A}_n$, $q_n$), for the candidate point $\x_{n+1}$ we get
\begin{equation*}
\begin{aligned}
q_{n+1}(\x_{n+1}, g_{n+1}) = \displaystyle\sum_{i=0}^{L} m_{n+1}(\x_{n+1}^q(B_i)) \mathbf{1}_{z \in B_i}.
\end{aligned}
\end{equation*}
\end{proposition}
Intuitively, the updated percentile is equal to the updated GP mean at one of the $\Xmc$ points, that depends on which interval $g_{n+1}$ (or equivalently, $z$) falls.

\begin{figure}[!ht]
\includegraphics[scale=0.3]{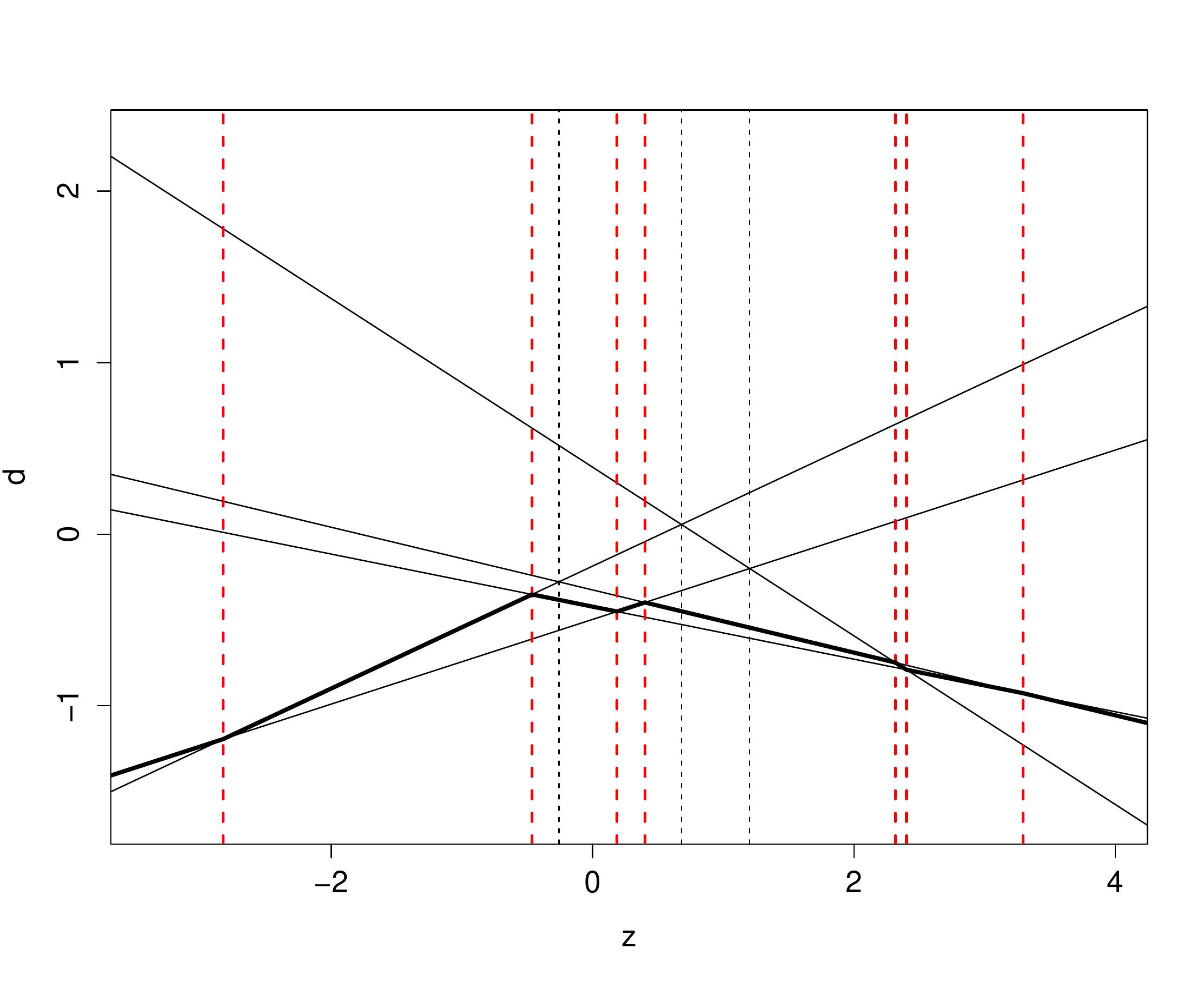}
\caption{Evolution of the percentile point as a function of the value of $z$. Each plain line represents a point of $\Xmc$, and the vertical lines the relevant intersections $I_i$. The second lowest line is shown in bold.}\label{fig:graphical}
\end{figure}

Figure \ref{fig:graphical} provides an illustrative example of this proposition for $l=5$, and $\alpha=40\%$. The values of $a$ and $b$ are given by a GP model, which allows us to draw the straight lines (black) as a function of $z$.
Each line corresponds to a point $\x_{\text{MC}}^i$. The intersections for which the percentile point changes are shown by the vertical lines.
For each interval, the segment corresponding to the percentile point (second lowest line) is shown in bold. We see that depending on the value of $z$ (that is, the value of $g_{n+1}$), the percentile point changes.
On the example, $j_i$ takes successively as values 2, 3, 1, 4, 3 and 5.

\begin{remark}
  Although the number of intersections grows quadratically with the MC sample size, 
  finding the set of percentile points can be done very efficiently, based on two important elements: first, 
  the number of distinct percentile points is much smaller than the number of intersections (five and ten, respectively, on Figure \ref{fig:graphical},
  but this difference increases rapidly with $l$); second, the order of the straight lines remains the same except for two elements for two adjacent intervals.
  This later feature allows us to avoid numerous calls to sorting functions.
  
  In the following, the notations $B_i$ and $L$ denote the effective intervals and number of intervals, respectively, that is, the intervals for which the percentile points are different. 
\end{remark}
%

\noindent {\bf 3.2. Infill criterion based on probability of exceedance}\label{sec:proba}

The proposition \ref{MAJpercentile} allows us to express the percentile estimator at step $n+1$ as a function of the candidate point $\x_{n+1}$ and corresponding value $g_{n+1}$.
In this section, we use this formulation to define a SUR criterion, that is, an uncertainty measure related to our estimator that can be minimized by a proper choice of $\x_{n+1}$.

This criterion is inspired from related work in probability of failure estimation \cite{bect} and multi-objective optimization \cite{victorMO}, that take advantage 
of the closed-form expressions of probabilities of exceeding thresholds in the GP framework.

By definition, the percentile is related to the probability of exceedance by
\begin{equation}
 \P(G(X) \geq q(G(X)))=1-\alpha.
\end{equation}

Our idea is the following. The probability $\mathbb{P}(G(\x) \geq q_n | \mathcal{A}_n)$, available for any $\x \in \Xset$, is 
in the ideal case ($G$ is exactly known) either zero or one, and, if $q_n = q(G(X))$, the proportion of ones is exactly equal to $1 - \alpha$.
At step $n$, a measure of error is then:
\begin{equation}
 J^{\text{prob}}_n= \left|\int_{\mathbb{X}} \P(G(\x) \geq q_n |\mathcal{A}_n) d\x -(1-\alpha) \right| = \left| \Gamma_n-(1- \alpha) \right|,
\end{equation}
with $\Gamma_{n}= \int_{\mathbb{X}} \P \left( G(\x) \geq q_n | \mathcal{A}_n \right) d\x $.
%
%
%

Following the SUR paradigm, the point we would want to add at step $n+1$ is the point $\x_{n+1}^*$ satisfying
\begin{equation}
 \x_{n+1}^* = \underset{\x_{n+1} \in \mathbb{X}}{\operatorname{argmin}\, }  J^{\text{prob}}_{n+1}(\x_{n+1}).
\end{equation}
As seen in proposition \ref{MAJpercentile}, $q_{n+1}$, and consequently $J^{\text{prob}}_{n+1}(\x_{n+1})$, depend on the candidate evaluation $g_{n+1}$, which makes it computable only by evaluating $g$. To circumvent this problem, we replace $g_{n+1}$ by its distribution conditional on $\A_n$.
%
We can then choose the following criterion to minimize (indexed by $\x_{n+1}$ to make the dependency explicit): 
\begin{equation}
  J^{\text{prob}}_{n}(\x_{n+1}) = \left| \E_{G_{n+1}} \left( \Gamma_{n+1}(x+1) \right) - (1- \alpha) \right| 
\end{equation}
where now,
\begin{equation}
\label{eq:gammaplus}
 \Gamma_{n+1}(x_{n+1})= \int_{\mathbb{X}} \P \left( G(\x) \geq q_{n+1} | A_{n+1} \right) d\x,
\end{equation}
with $A_{n+1}= \mathcal{A}_{n+1} \cup \left(\x_{n+1}, G_{n+1} \right)$ and $G_{n+1}$ is still in its random form. 

We show that
\begin{proposition}
\label{proba}
Using previous notations and under our first strategy, 
\begin{equation*}
\begin{aligned}
&J^{\text{prob}}_{n}(\x_{n+1}) = \\
&\Bigg| \int_{\mathbb{X}} \displaystyle\sum_{i=1}^{L-1} \Bigg[ \Phi_{r_i^n} \left(e^i_n(\x_{n+1}; \x) , f_n^i(\x_{n+1}, I_{i+1})\right) - \Phi_{r_i^n(\x_{n+1}, \x)} \left(e^i_n(\x_{n+1}; \x) , f_n(\x_{n+1}, I_{i})\right) \\
& + \Phi_{r_i^n}\left((e^i_n(\x_{n+1}; \x) , f_n^i(\x_{n+1}, I_{1})\right) + \Phi_{-r_i^n}\left(e^i_n(\x_{n+1}; \x) , -f_n^i\left((\x_{n+1}, I_{L}\right) \right) \Bigg] d\x - (1- \alpha) \Bigg| \\
\end{aligned}
\end{equation*}

where

$$e_n^i(\x_{n+1}; \x; \x_{n+1}^q(B_i))= \frac{m_n(\x) - m_n({\x_{n+1}^q(B_i)})}{\sigma_{W}}, \, \, f_n^i(\x_{n+1}; I_i)= I_is_n(\x_{n+1}),$$

$$\sigma_{W}=s_n(\x)^2+\frac{k_n(\x^q_{n+1}(B_i), \x_{n+1})^2}{s_n(\x_{n+1})^2} - 2 \frac{k_n(\x^q_{n+1}(B_i), \x_{n+1})k_n(\x, \x_{n+1})}{s_n(\x_{n+1})^2}$$

and $\Phi_{r^i_n}$ is the cumulative distribution function (CDF) of the centered Gaussian law of covariance matrix 

$$K=  \begin{pmatrix} 1 & r^i_n \\
r^i_n & 1 \\
\end{pmatrix}$$

with $$r^i_n=\frac{k_n(\x_{n+1}^q(B_i), \x_{n+1})- k_n(\x, \x_{n+1})}{\sqrt{s_n(\x)^2 + \frac{k_n(\x_{n+1}^q(B_i), \x_{n+1})^2}{s_n(\x_{n+1})^2} - 2 \frac{k_n(\x_{n+1}^q(B_i), \x_{n+1})k_n(\x, \x_{n+1})}{s_n(\x_{n+1})^2}} s_n(\x_{n+1})}.$$

\end{proposition}

The proof is deferred to the Appendix.

Despite its apparent complexity, this criterion takes a favourable form, since it writes as a function of GP quantities at step $n$
($m_n$, $s_n$ and $k_n$), which can be computed very quickly once the model is established. Besides, it does not require
conditional simulations (as the criterion in \cite{jala2}), which is a decisive advantage both in terms of 
computational cost and evaluation precision.

Let us stress here, however, that evaluating this criterion requires a substantial computational effort, as it 
takes the form of an integral over $\Xset$, which must be done numerically. An obvious choice here is to use 
the set $\Xmc$ as integration points. 
Also, it relies on the bivariate Gaussian CDF, which also must be computed numerically. Very efficient programs can be found, 
such as the R package \texttt{pbivnorm} \cite{kenkel2011pbivnorm}, which makes this task relatively inexpensive. 

\noindent {\bf 3.3. Infill criterion based on the percentile variance}\label{sec:variance}

Accounting for the fact that, although not relying on conditional simulations,
$J^{\text{prob}}$ is still expensive to compute, 
we propose here an alternative, that does not require numerical integration over $\Xset$.



Since we want $q_n$ to converge to the quantile, it is important that this estimator becomes increasingly \textit{stable}. 
The variance of $q_{n+1}|A_{n+1}$ is a good indicator of this stability, as it gives the fluctuation range of $q_{n+1}$ 
as a function of the different possible values of $g_{n+1}$. 
However, choosing the point that minimizes at each step $n$ this variance has no sense here, as choosing
$\x_{n+1} \in \{\x_1, \ldots, \x_n\}$ (that is, duplicating an existing observation) would result in $\Var(q_{n+1}|A_{n+1})=0$.

Inversing the SUR paradigm, we propose to choose the point that \emph{maximizes} this variance.
By doing so, we will obtain the sequence of points which evaluations have a large impact on the estimator value,
hence reducing sequentially the instability of our estimator:
%
%
%
%
\begin{equation}\label{eq:jvar}
 J^{\text{Var}}_{n}(\x_{n+1}) = \Var_{G_{n+1}}( q_{n+1} | A_{n+1})
\end{equation}
where once again $A_{n+1}$ denotes the conditioning on $\A_n \cup (\x_{n+1}, G_{n+1})$, with $G_{n+1}$ random. 
We can show that:
\begin{proposition}\label{var}
Using the previous notations, conditionally on $\A_n$ and on the choice of $\x_{n+1}$:
\begin{multline}
J^{\text{Var}}_{n}(\x_{n+1}) =
\displaystyle\sum_{i=1}^L \left[ k_n(\x^q_{n+1}(B_i), \x_{n+1})\right]^2 V(s_n(\x_{n+1}), I_{i+1}, I_i) P_i \nonumber \\
 + \displaystyle\sum_{i=1}^L \left[ m_n(\x_{n+1}^q(B_i) - k_n(\x^q_{n+1}(B_i), \x_{n+1})E(s_n(\x_{n+1}), I_{i+1}, I_i)\right]^2 \left( 1 - P_i \right) P_i\nonumber \\
- 2 \displaystyle\sum_{i=2}^L \displaystyle\sum_{j=1}^{i-1} \left[m_n(\x_{n+1}^q(B_i) - k_n(\x^q_{n+1}(B_i), \x_{n+1})E(s_n(\x_{n+1}), I_{i+1}, I_i)\right] P_i \nonumber\\
\left[m_n(\x_{n+1}^q(B_i) - k_n(\x^q_{n+1}(B_i), \x_{n+1})E(s_n(\x_{n+1}), I_{j+1}, I_j)\right] P_j,
\end{multline}
where:
$$P_i=\Phi(s_n(\x_{n+1})I_{i+1})-\Phi(s_n(\x_{n+1})I_i),$$
$$E(s_n(\x_{n+1}), I_{i+1}, I_i)= \frac{1}{s_n(\x_{n+1})} \left( \frac{\phi(s_n(\x_{n+1})I_{i+1})-\phi(s_n(\x_{n+1})I_i)}{\Phi(s_n(\x_{n+1})I_{i+1})- \Phi(s_n(\x_{n+1})I_i)} \right),$$
and
\begin{multline*}
 V(s_n(\x_{n+1}), I_{i+1}, I_i)= \\
 \frac{1}{s_n(\x_{n+1})^2} \left[ 1+ \frac{s_n(\x_{n+1}) \phi( I_{i+1}) - s_n(\x_{n+1}) \phi( I_{i})}{\Phi( I_{i+1})-\Phi( I_{i})} - \left( \frac{\phi(I_{i+1}) - \phi(I_i)}{\Phi(I_{i+1}) - \Phi(I_i)} \right)^2\right],
\end{multline*}
%
     
for $\Phi$ and $\phi$ respectively the cumulative distribution function and density function of the standard Gaussian law. 

\end{proposition}
The proof is deferred to Appendix.

Again, this criterion writes only as a function of GP quantities at step $n$ ($m_n$, $s_n$ and $k_n$).
As it does not require numerical integration nor the bivariate CDF, it is considerably cheaper
to compute than the previous one.

\begin{figure}[!ht]
\includegraphics[width=\textwidth]{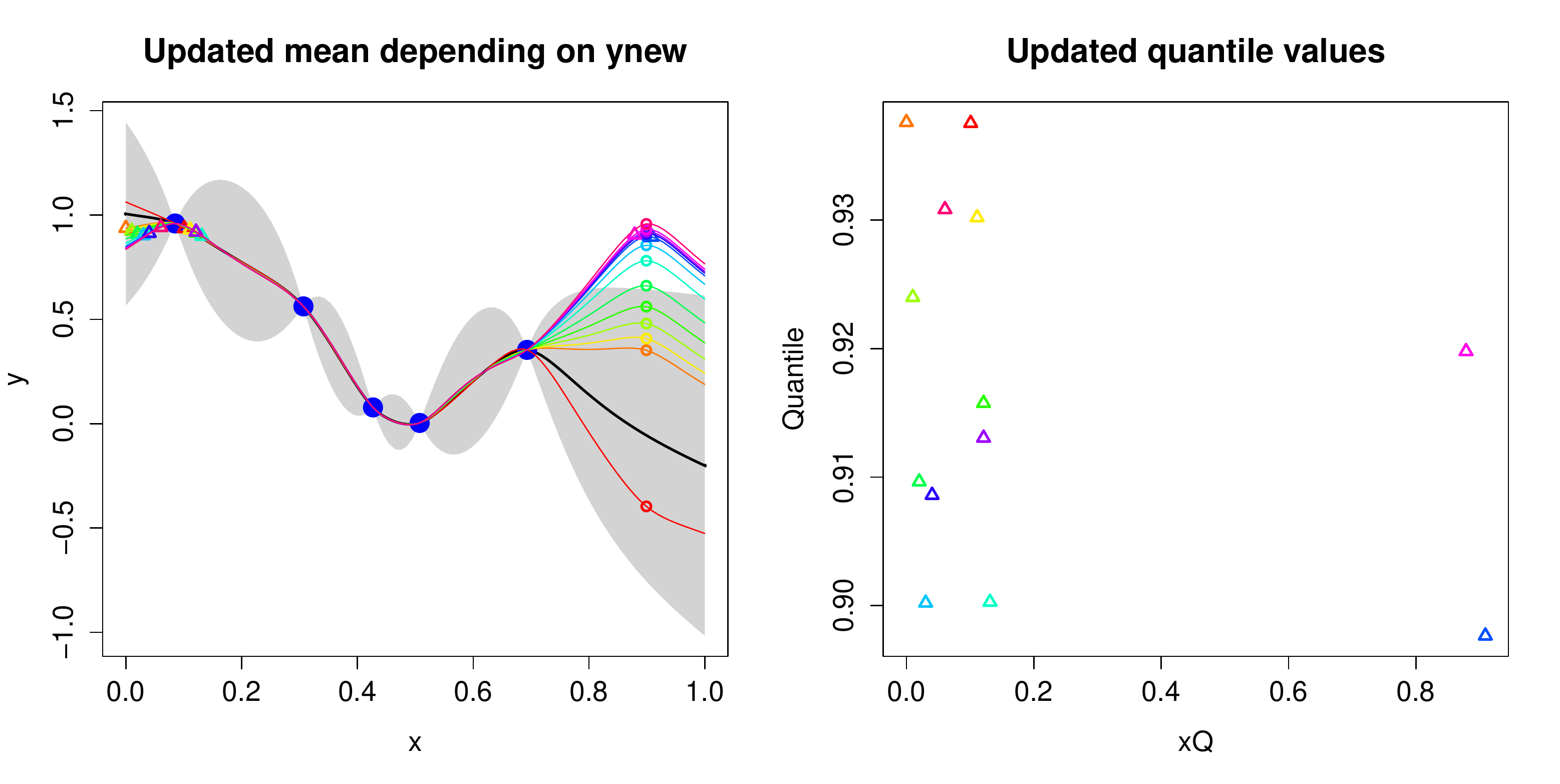}
\caption{Illustration of the $J^{\text{Var}}$ criterion. Left: GP model (black bold line and grey area) and updated GP mean (other lines) depending on the value of $g_{new}$ (circles) for $\x_{n+1}=0.9$. 
The corresponding 90\% percentiles $q_{n+1}$ are shown with the triangles. Right: percentile values only, indexed by the corresponding percentile points.}\label{fig:illusvarQ}
\end{figure}

Figure \ref{fig:illusvarQ} provides an illustration of the concepts behind this criterion, by showing how 
different values of $g_{n+1}$ affect the estimator. Here, one updated mean is drawn for each interval $B_i$
(that is, with $z$ taking its value in the middle of the interval). The corresponding 90\% percentiles, 
as well as the percentile points $\x_{n+1}^q$ vary substantially, depending on $g_{n+1}$, which results in a 
large variance $\Var(q_{n+1}|A_{n+1})$. Hence, the point $\x_{n+1}=0.9$ can be considered as highly informative for our estimator.

\section{Numerical simulations}\label{sec:exp}

\noindent {\bf 4.1. Two-dimensional example}\label{sec:2dim}

As an illustrating example, we use here the classical Branin test function (see \cite{dixon1978towards} Equation \ref{eq:branin} in Appendix).
On $[0,1]^2$, the range of this function is approximately $[0, 305]$.

We take: $X_1, X_2 \sim \mathcal{U}[0,1]$, and search for the $85\%$ percentile. The initial set of experiments consists of 
seven observations generated using Latin Hypercube Sampling (LHS), and 11 observations are added sequentially using both SUR
strategies. The GP models learning, prediction and update is performed using the \texttt{R} package \texttt{DiceKriging} \cite{roustant2012dicekriging}.
The covariance is chosen as Mat\'ern $3/2$ and the mean as a linear trend.

For $\Xmc$, we used a 1000-point uniform sample on $[0,1]^2$.
For simplicity purpose, the search of $\x_{n+1}$ is performed on $\Xmc$, although a continuous optimizer algorithm could have been used here.
The actual percentile is computed using a $10^5$-point sample.

\begin{figure}[!ht]
\includegraphics[trim=0 98mm 0 0mm, clip, width=.66\textwidth]{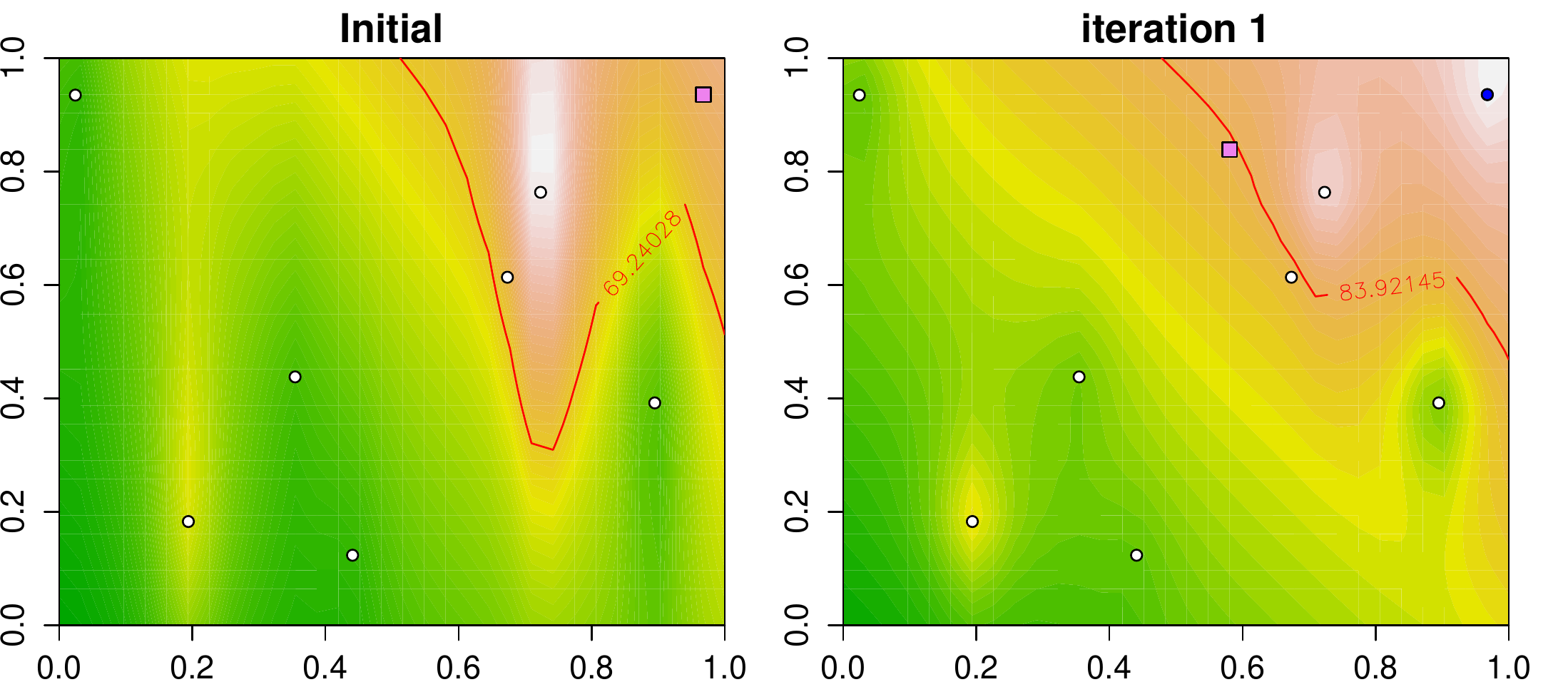}
\includegraphics[trim=0mm 0mm 112mm 98mm, clip, width=.33\textwidth]{2D_7plus11_var_4ite.pdf}
\caption{Contour lines of the GP mean and experimental set at three different $n$ values (7, 8, and 15) with the $J_n^{\text{Var}}$ criterion. 
The initial observations are shown with white circles, 
the observations added by the sequential strategy with blue circles, and the next point to evaluate with violet squares.
The line shows the contour corresponding to the percentile estimate.}\label{fig:doe}
\end{figure}

\begin{figure}[!ht]
 \includegraphics[trim=0 0 0 10mm, clip, width=.49\textwidth]{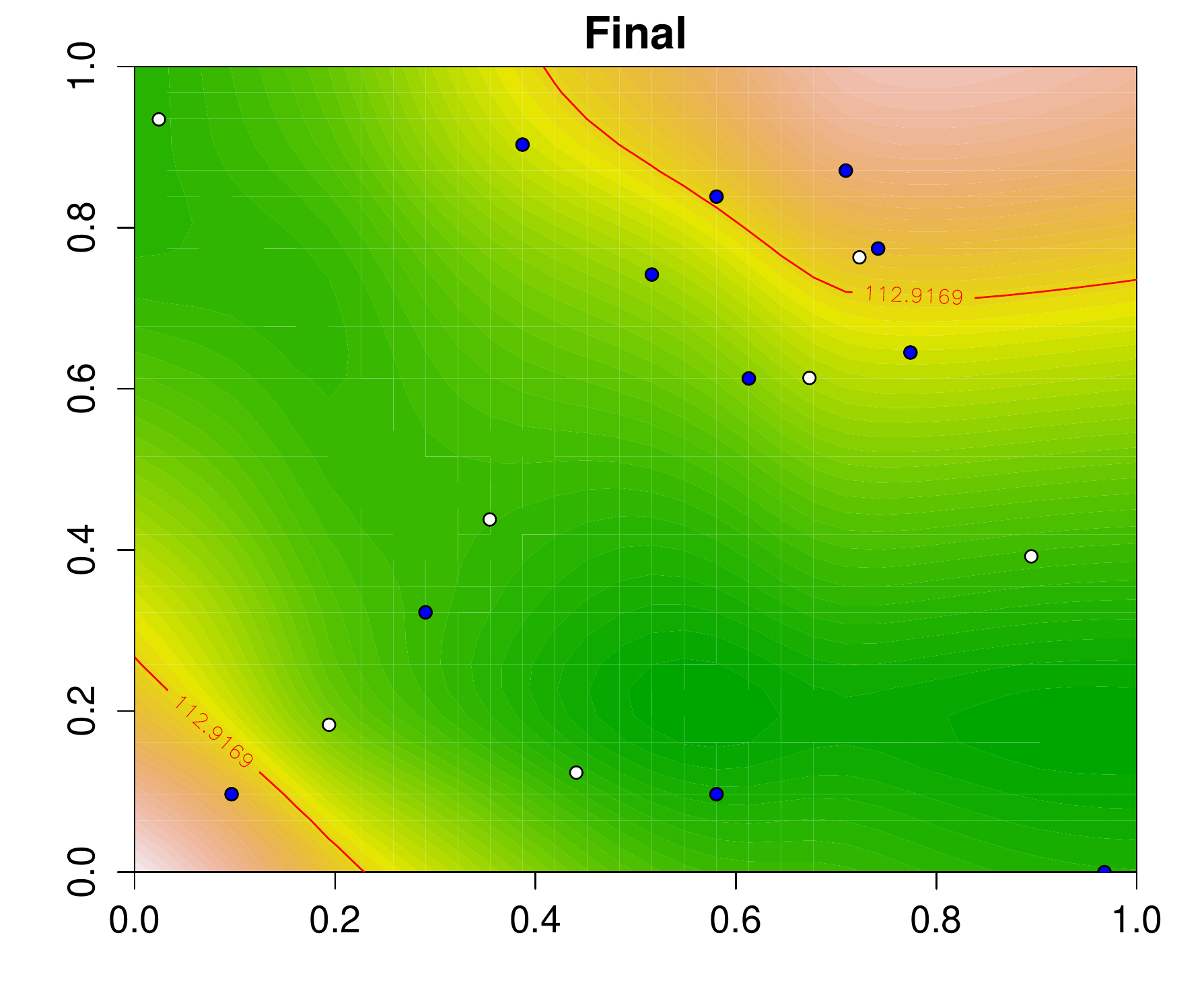}
 \includegraphics[trim=0 0 0 10mm, clip, width=.49\textwidth]{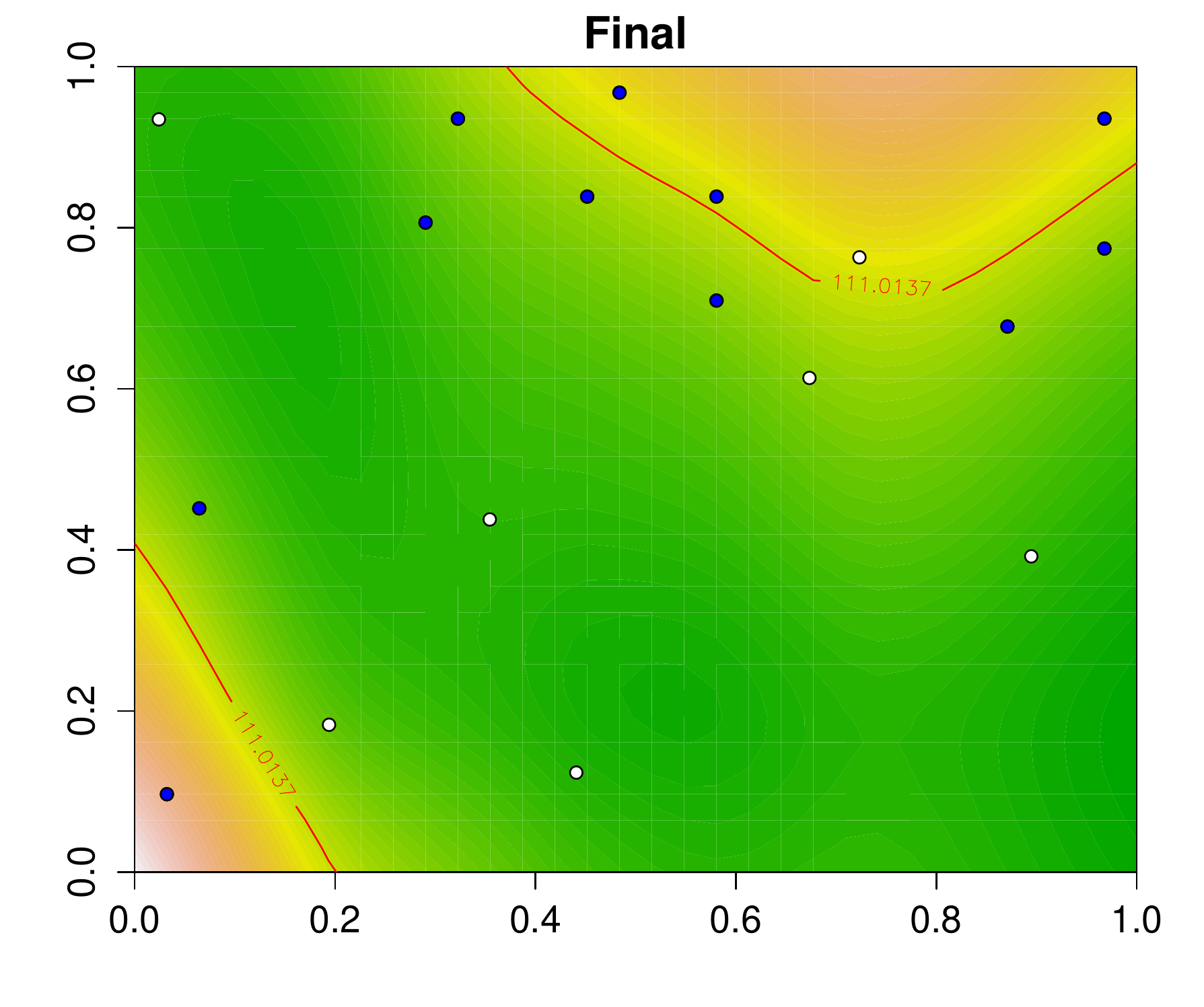}
   \caption{Comparison of observation sets obtained using $J_n^{\text{prob}}$ (left) and $J_n^{\text{Var}}$ (right).}\label{fig:2doe}
\end{figure}

\begin{figure}[!ht]
\includegraphics[width=.49\textwidth]{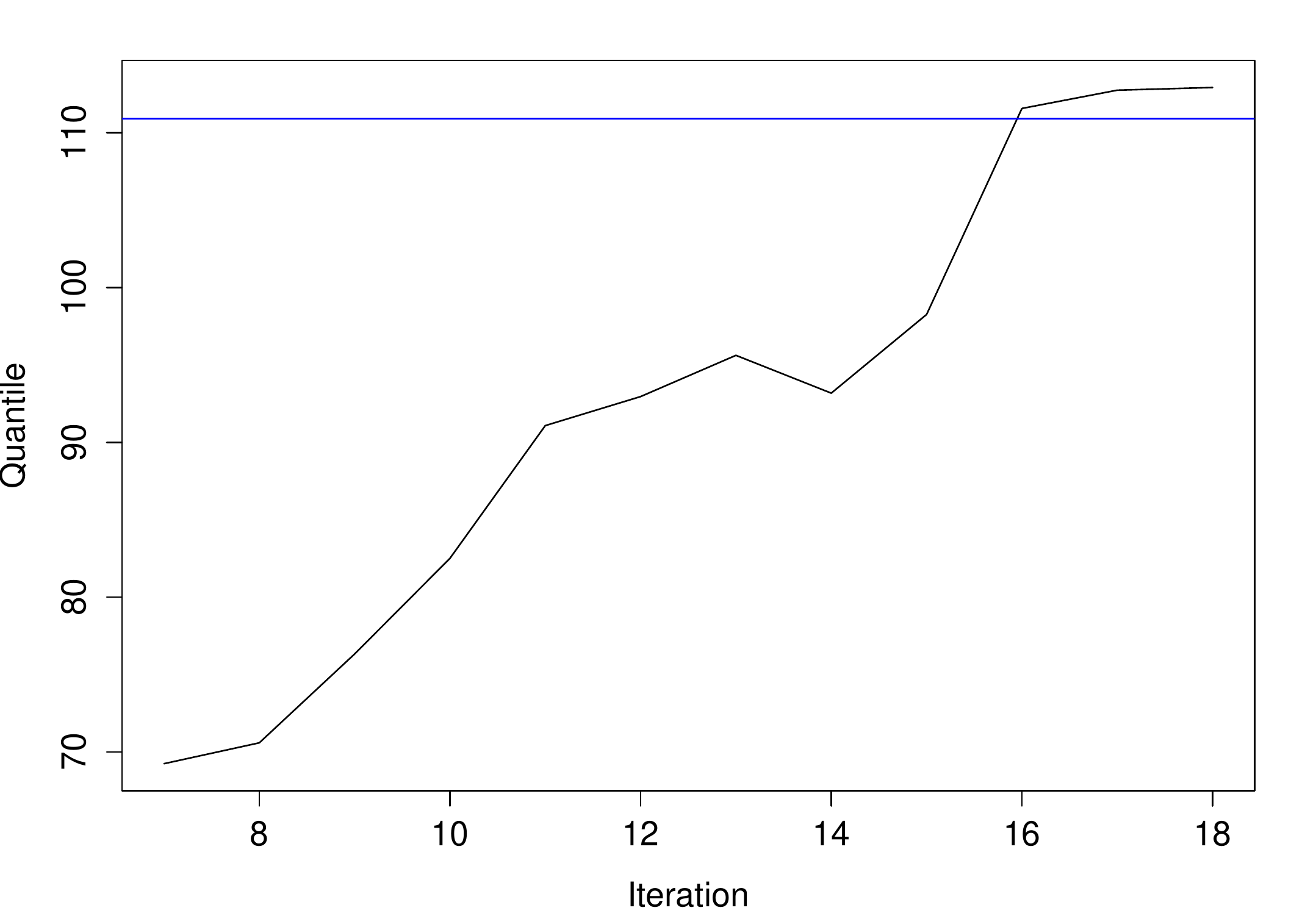}
\includegraphics[width=.49\textwidth]{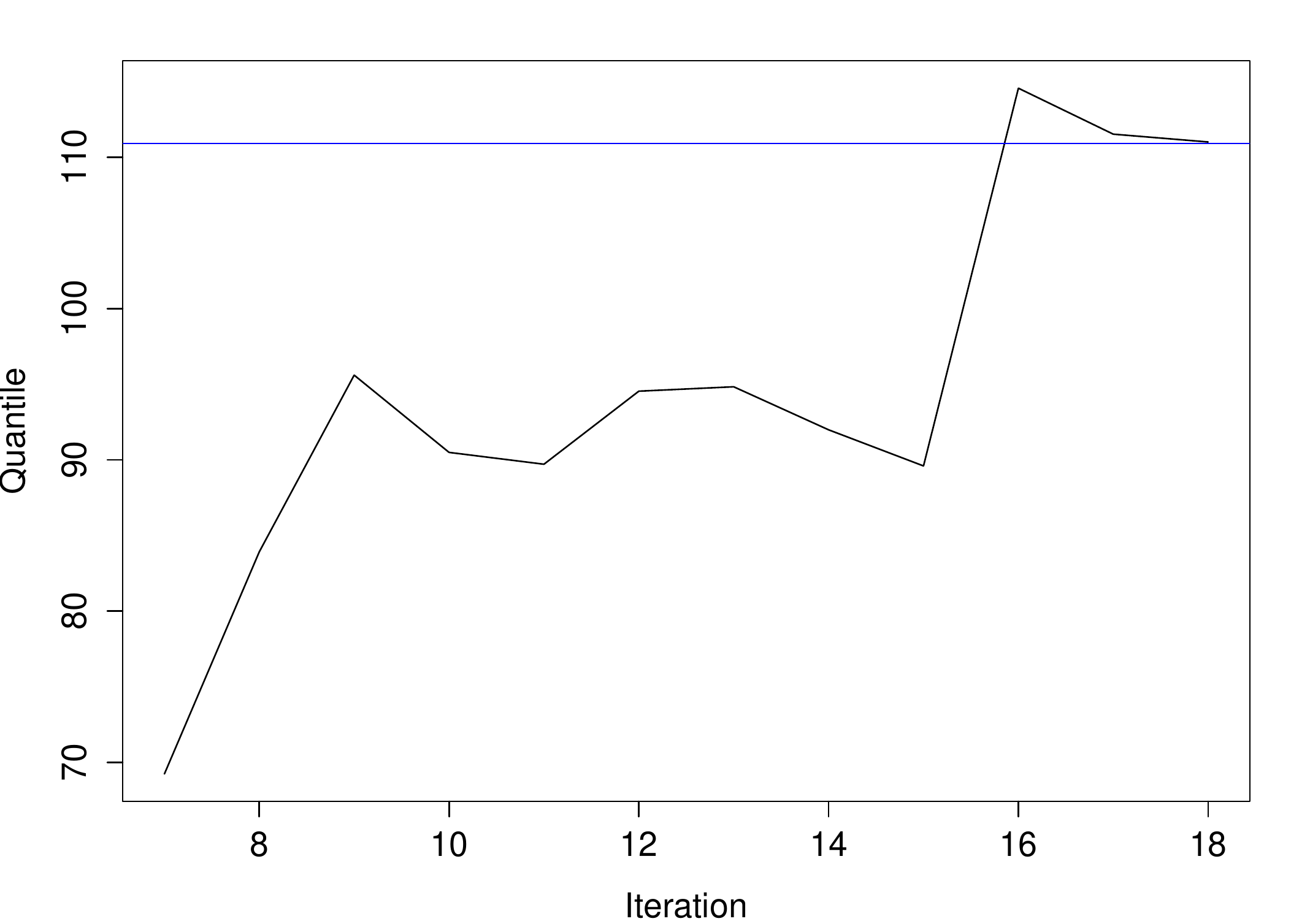}
\caption{Evolution of the percentile estimates using $J_n^{\text{prob}}$ (left) and $J_n^{\text{Var}}$ (right) for the 2D problem. The horizontal line shows the actual $85^{th}$ percentile.}\label{fig:quantileEvol2D}
\end{figure}

Figure \ref{fig:2doe} reports the final set of experiments, along with contour lines of the GP model mean,
and Figure \ref{fig:quantileEvol2D} the evolution of the estimators. In addition, Figure \ref{fig:doe} 
shows three intermediate stages of the $J_n^{\text{Var}}$ run.

Figure \ref{fig:doe} reveals the dynamic of our strategy: from the initial design of experiments, the top right corner of the domain
is identified as the region containing the highest 15\% values. Several observations are added in that region until the kriging 
approximation becomes accurate, then a new region (bottom left corner) is explored (square point, Figure \ref{fig:doe} right).

The two strategies lead to relatively similar observation sets (Figure \ref{fig:2doe}), that mostly consist of values close to the contour line 
corresponding to the $85^{th}$ percentile (exploitation points), and a couple of space-filling points (exploration points).
With 18 observations, both estimators are close to the actual value (in particular with respect to the range of the function),
yet additional observations may be required to achieve convergence (Figure \ref{fig:quantileEvol2D}).

\noindent {\bf 4.2. Four and six dimensional examples}\label{sec:4dim}

We consider now two more difficult test functions, with four and six dimensions, respectively
(\textit{hartman} and \textit{ackley} functions, see Equations \ref{eq:hartman} and \ref{eq:ackley} in Appendix).
Both are widely used to test optimization strategies \cite{dixon1978towards}, and are 
bowl-shaped, multi-modal functions.

We take on both cases: $\X \sim \Nlaw \left( \frac{1}{2}, \Sigma \right)$, with $\Sigma$ a symmetric matrix
with diagonal elements equal to $0.1$ and other elements equal to $0.05$.
The initial set of observations is taken as a 30-point LHS,
and 60 observations are added sequentially. A 3000-point sample from the distribution of $\X$ is used for $\Xmc$
(renewed at each iteration),
and the actual percentile is computed using a $10^5$-point sample.
Again, the GP covariance is chosen as Mat\'ern $3/2$ and the mean as a linear trend.

The criteria are optimized as follow: a (large) set of $10^5$ candidates is generated from the distribution of $\X$,
out of which a shorter set of $300$ ``promising'' points is extracted. Those points are drawn randomly from the large
set with weights equal to $\phi \left( \frac{q_n - m_n(\x)}{s_n(\x)} \right)$. Hence, higher weights are given to points
either close to the current estimate and/or with high uncertainty. The criterion is evaluated on this 
subset of points and the best is chosen as the next infill point.
In addition, for $J_n^{\text{Var}}$ a local optimization is performed, starting from the best point of the subset
(using the BFGS algorithm, see \cite{algo}). Due to computational constraints, this step is not applied to $J_n^{\text{prob}}$,
which is more costly. However, preliminary experiments have shown that only a limited gain 
is achieved by this step.

As an baseline strategy for comparison purpose, we also include a ``random search'', that is,
the $\x_{n+1}$ are sampled randomly from the distribution of $X$.

Several percentile levels are considered in order to cover a variety of situations: $5\%$ and $97\%$ for the 4D
problem and $15\%$ and $97\%$ for the 6D problem. Due to the bowl-shape of the functions, low levels are 
defined by small regions close to the center of the support of $X$, while high levels correspond to the 
edges of the support of $X$. Besides, it is reasonable to assume that levels farther away from $50\%$ are 
more difficult to estimate.

As an error metric $\varepsilon$, we consider the absolute difference between the percentile estimator and its actual value.
We show this error as a percentage of the variation range of the test function. Since $X$ is not bounded, 
the range is defined as the difference between the $0.5$ and $99.5$ percentiles of $g(X)$.

\begin{figure}[!ht]
\includegraphics[width=.49\textwidth]{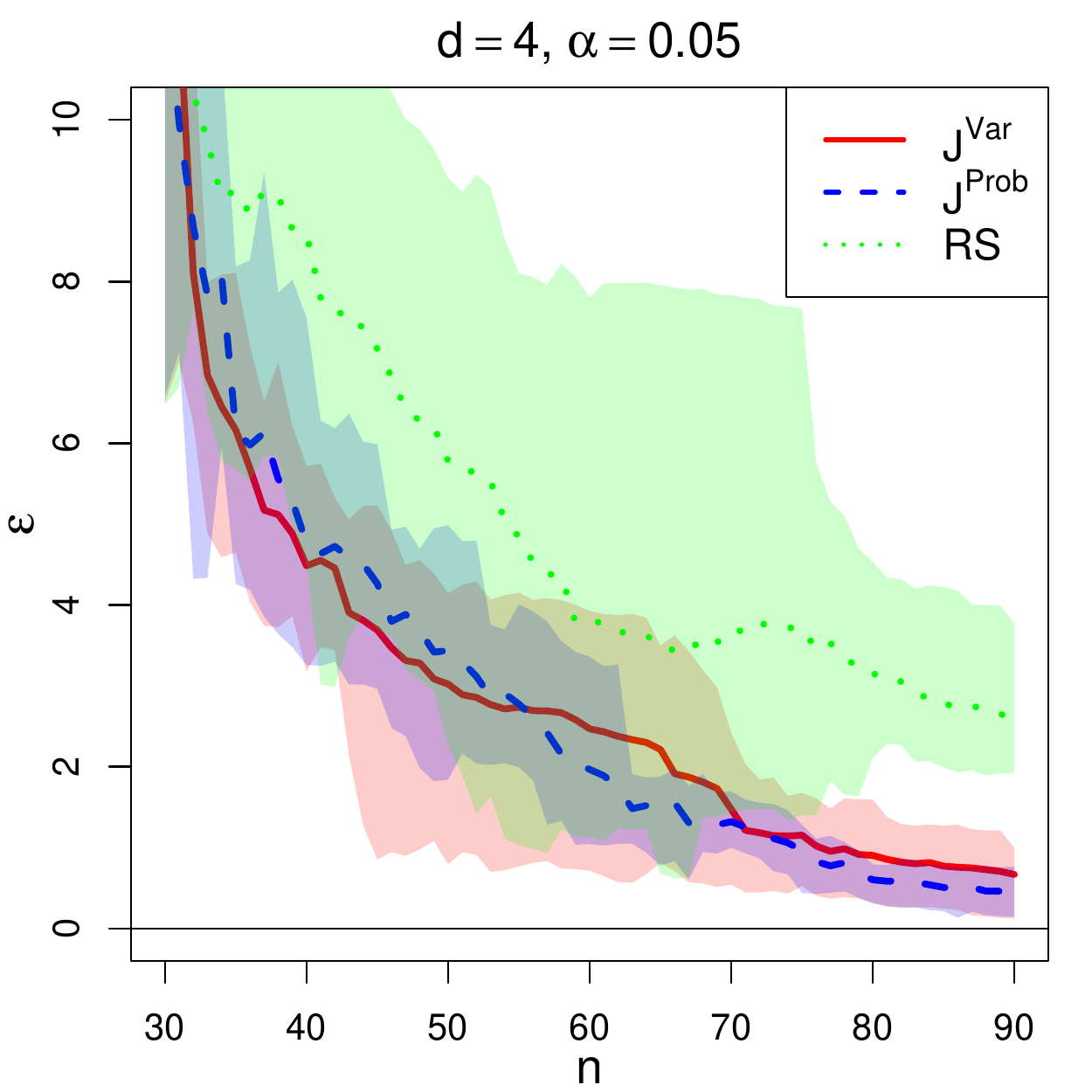}
\includegraphics[width=.49\textwidth]{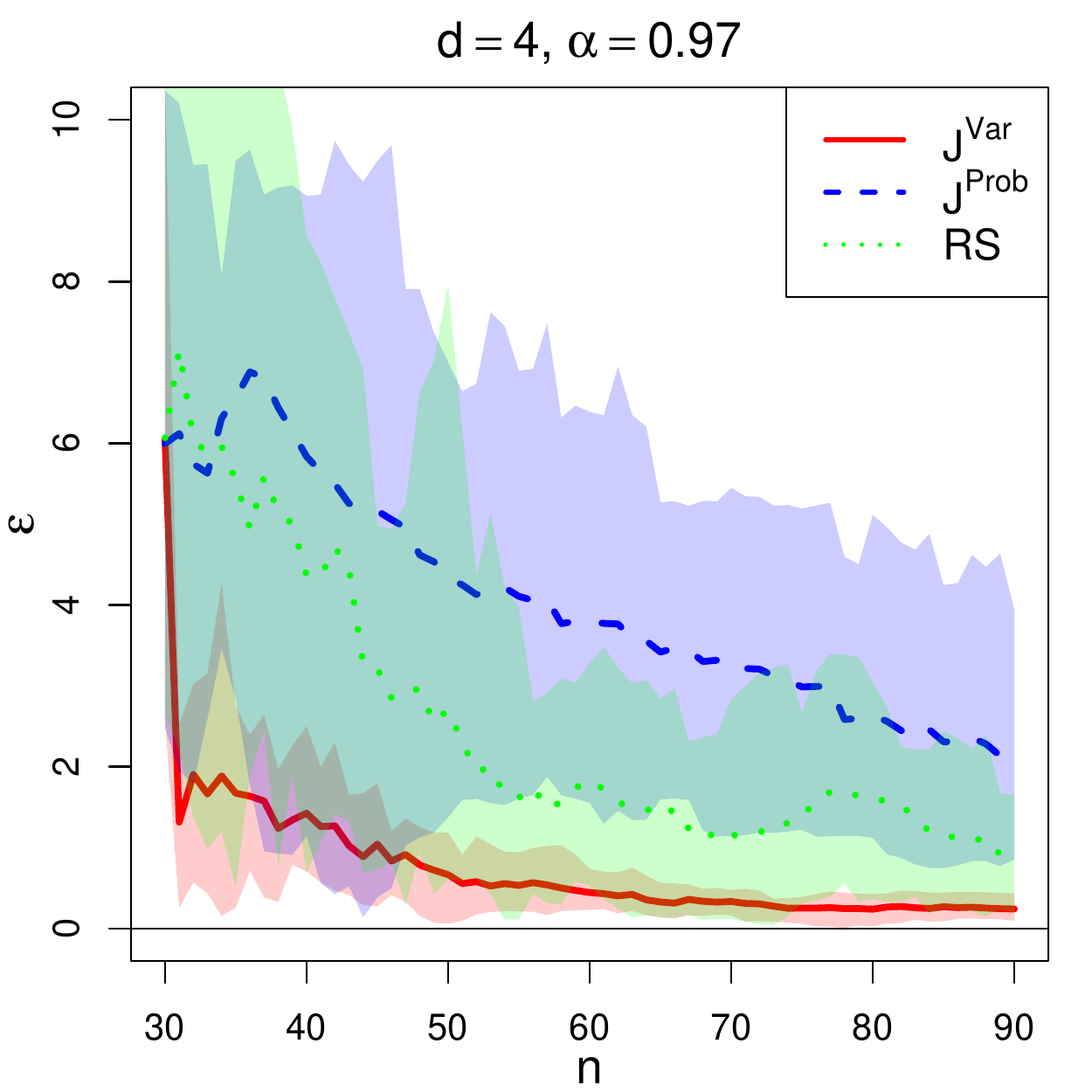}
\includegraphics[width=.49\textwidth]{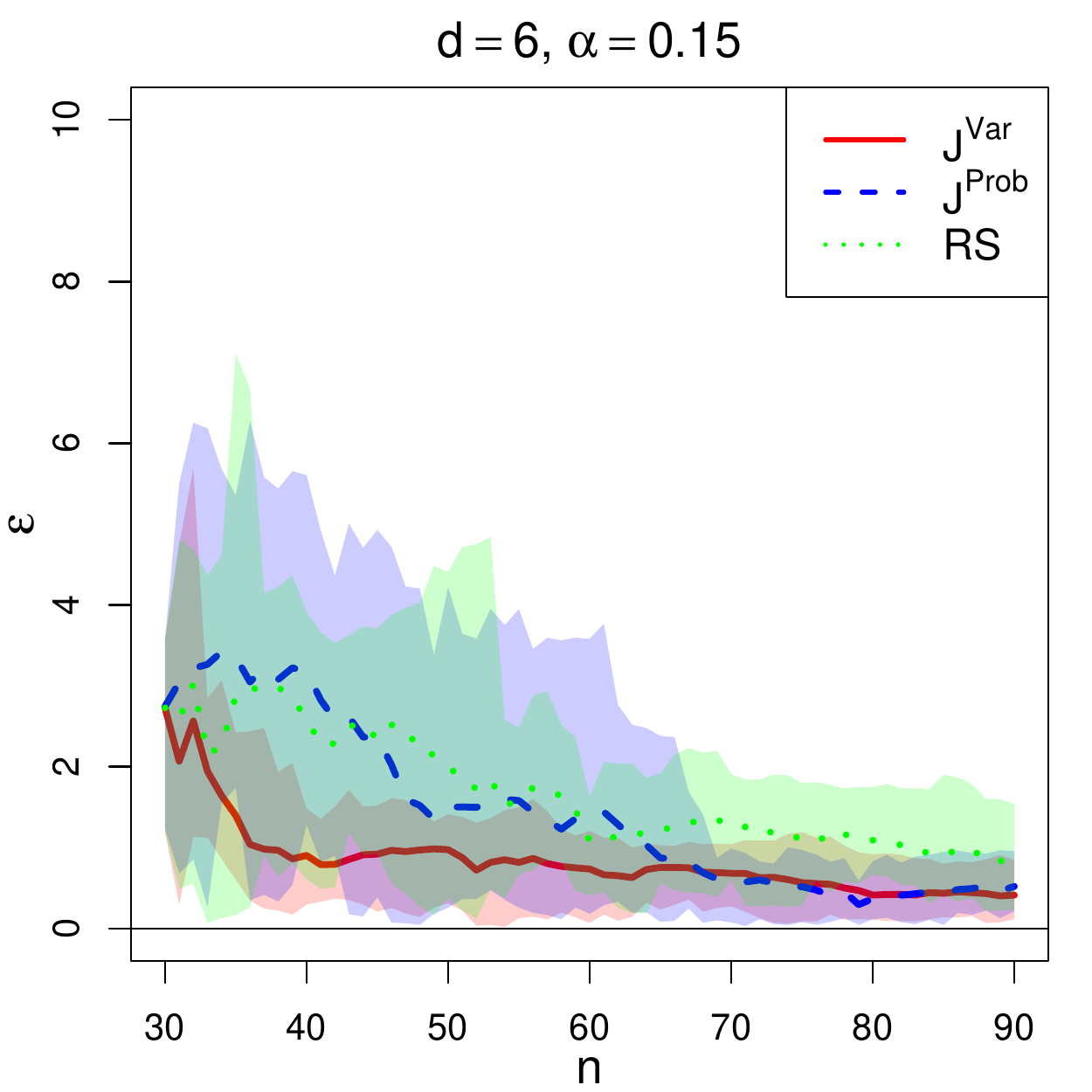}
\includegraphics[width=.49\textwidth]{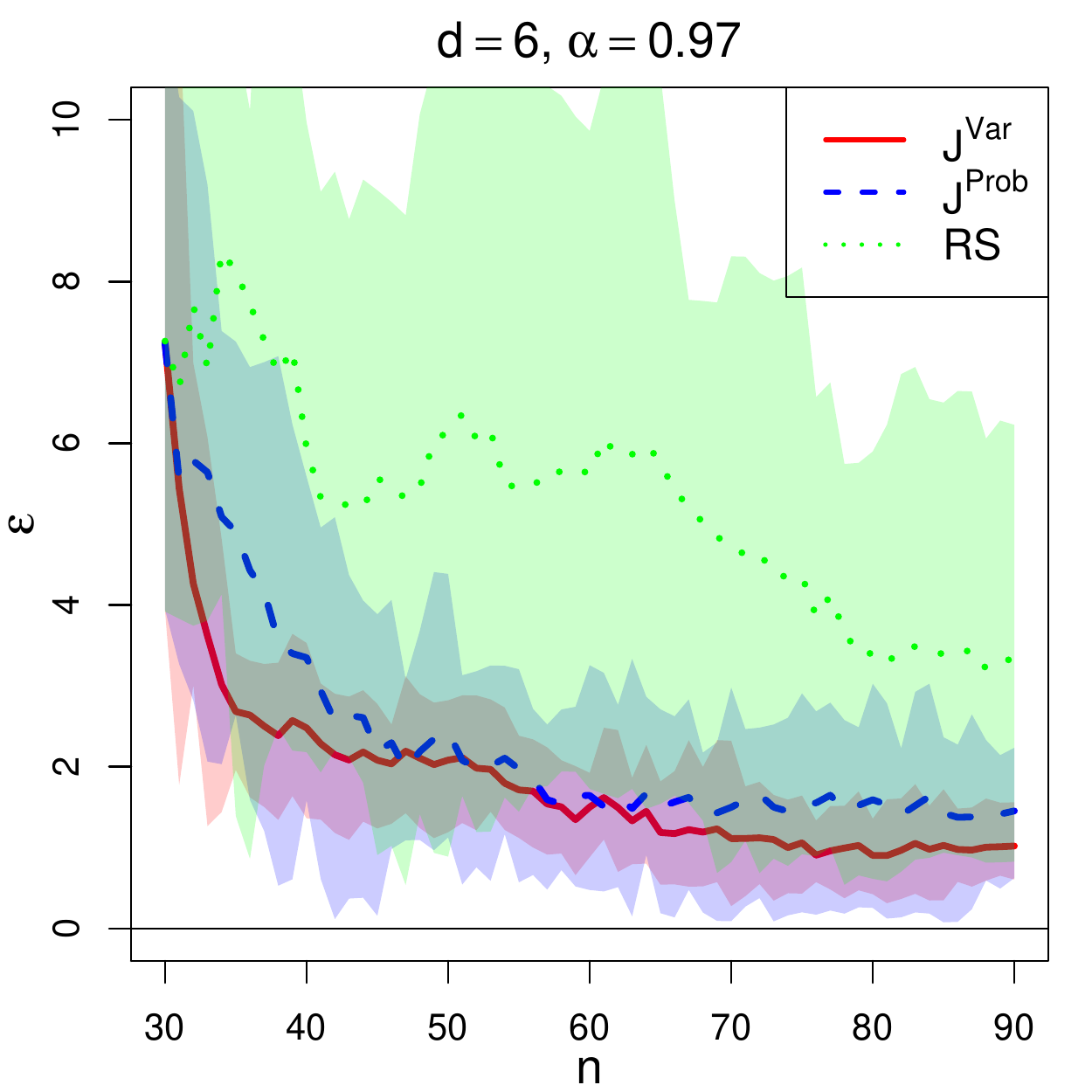}
\caption{Evolution of the percentile estimates using $J_n^{\text{prob}}$ (dashed line), $J_n^{\text{Var}}$ 
(plain line) or random search (RS, dotted line) for the 4D and 6D problems and several percentile levels. 
The lines show the average error and the shaded areas the 10\% and 90\% quantile errors over the runs.}\label{fig:quantileEvol4D}
\end{figure}

To assess the robustness of our approach, the experiments are run ten times for each case, starting with a different initial 
set of observations. The evolution of the estimators (average, lowest and highest error metric values over the ten runs) 
is given in Figures \ref{fig:quantileEvol4D}. 

First of all, we see that except on one case (4D, $\alpha=0.97$ and $J_n^{\text{prob}}$),
on average both strategies provide estimates with less than $2\%$ error after approximately 
30 iterations (for a total of 60 function evaluations), which plainly justifies the use 
of GP models and sequential strategies in a constrained budget context.

For $d=4$, $\alpha=0.05$, both methods seem to converge to the actual percentile, 
 $J_n^{\text{prob}}$ being slightly better, in particular in terms of consistency
 and for the latest steps.
 
For $d=4$, $\alpha=0.97$, $J_n^{\text{Var}}$ reaches very quickly for all runs
a good estimate (less that $1\%$ error), yet seems to converge then slowly to the exact solution.
This might be explained by the relative mismatch between the GP model and the test function.
$J_n^{\text{prob}}$ performs surprisingly poorly; we conjecture that a more exploratory behavior
compared to $J_n^{\text{Var}}$ hinders its performance here.

For $d=6$, $\alpha=0.15$, both approaches reach consistently less than $1\%$ error. However,
they outperform only moderately the random search strategy here. This might indicate that
for central percentile values, less gain can be achieved by sequential strategies, as a large
region of the design space needs to be learned to characterize the percentile, making 
space-filling strategies for instance competitive.

Finally, for $d=6$, $\alpha=0.97$, both approaches largely outperform random search,
yet after a first few very efficient steps seem to converge only slowly to the actual percentile.

In general, those experiments show the ability of our approach to handle multi-modal 
black-box functions, with input space dimensions typical of GP-based approaches. 
Our results seem to indicate a better efficiency of the $J_n^{\text{Var}}$ criterion,
yet the better convergence with $J_n^{\text{prob}}$ for $d=4$, $\alpha=0.05$ might
call for hybrid strategies, with early steps performed with $J_n^{\text{Var}}$
and late steps with $J_n^{\text{prob}}$.



%

\section{Concluding comments}\label{sec:conclusion}

We have proposed two efficient sequential Bayesian strategies for percentile estimation.
Both approaches rely on the analytical update formula for the GP-based estimator, 
which has been obtained thanks to the particular form of the GP equations and the
introduction of the \emph{quantile point} concept. Two criteria have then been proposed based either 
on probability of exceedance or on variance, for which closed-form expression have
been derived, hence avoiding the use of computationally intensive conditional simulations.
Numerical experiments in dimensions two to six have demonstrated the potential of both approaches.

There are of course some limitations of the proposed method, that call for future improvements.
Both strategies rely on the set $\Xmc$, which size is in practice limited by the computational resources
to a couple of thousands at most. This may hinders the use of our method for extreme percentile estimation, 
or for highly multi-modal functions. Combining adaptive sampling strategies or subset selection methods 
with our approaches may prove useful in this context.

Accounting for the GP model error (due to an inaccurate estimation of its hyper-parameters or a poor choice 
of kernel) is also an important task, that may improve greatly the efficiency and robustness of the approach.
Embrassing a fully Bayesian approach (as for instance in \cite{kennedy2001bayesian,gramacy2008bayesian})
may help address this issue, yet at the price of additional computational expense.

\section*{Aknowledgements}
The authors would like to thank Damien Leroux for his help on the efficient 
implementation of the algorithms on R.

\noindent {\bf 6.1. Proof of Proposition \ref{proba}}\label{sec:proofproba}

In the following, we denote $\E_n$ and $\P_n$ the expectation and the probability conditionally on the event $\mathcal{A}_n$.
Starting from Equation \ref{eq:gammaplus}, we have:
\begin{equation*}
\begin{aligned}
\E_{G_{n+1}}\left(\Gamma_{n+1}(\x_{n+1})\right) & = \E_{G_{n+1}} \left[  \int_{\mathbb{X}} \P(G(x) \geq q_{n+1})  | A_{n+1})dx \right] \\
& = \int_{\mathbb{X}}\E_{G_{n+1}} \left[ \E_{n} \left[ \mathbf{1}_{G(x) \geq q_{n+1}(x_{n+1})} | G_{n+1} \right] \right] dx  \\
&= \int_{\mathbb{X}} \E_n \left[ \mathbf{1}_{G(x) \geq q_{n+1}(\x_{n+1})} \right] dx \\
& = \int_{\mathbb{X}}\P_n( G(x) \geq q_{n+1}(\x_{n+1})) dx \;.
\end{aligned}
\end{equation*}
We get then:
\begin{equation*}
\begin{aligned}
J^{\text{prob}}_{n+1}(\x_{n+1})&= \left| \int_{\mathbb{X}} \P_n( G(x) \geq q_{n+1}(\x_{n+1}))dx -(1- \alpha) \right|\;.\\
\end{aligned}
\end{equation*}

Now, to get an closed form of our criterion, we have to develop $\P_n( G(x) \geq q_{n+1}(x_{n+1}))$.
To do so, we use Proposition \ref{MAJpercentile}. Denoting $Z= \frac{G_{n+1}-m_n(\x_{n+1})}{s_n(\x_{n+1})^2}$, we have:
\begin{equation*}
\begin{aligned}
\E_n  \left( \mathbf{1}_{G(x) \geq q_{n+1}(\x_{n+1})} \right) & = \displaystyle\sum_{i=0}^{L} \E_n \left[ \mathbf{1}_{G(x) \geq m_{n+1}(\x^q_{n+1}(B_i)} \mathbf{1}_{Z \in B_i} \right] \\
&=\displaystyle\sum_{i=1}^{L-1} \left( \P_n \left[ G(x) \geq m_{n+1}(\x^q_{n+1}(B_i)) \cap Z \leq I_{i+1} \right]\right. \\
& -\left. \P_n \left[ G(x) \geq m_{n+1}(\x^q_{n+1}(B_i)) \cap Z \leq I_{i} \right) \right] ) \\
& + \P_n(G(x) \geq m_{n+1}(\x^q_{n+1}(B_1)) \cap Z \leq I_1) \\
& + \P_n(G(x) \geq m_{n+1}(\x^q_{n+1}(B_L)) \cap Z \geq I_{L})\;. \\
\end{aligned}
\end{equation*} 

Now,
\begin{equation*}
\begin{aligned}
T_n:&= \P_n \left( G(x) \geq m_{n+1}(\x_{n+1}^q(B_i)) \cap  Z \leq I_i \right) \\
& = \P_n \left( m_{n+1}(\x_{n+1}^q(B_i))-G(x) \leq 0 \cap Z \leq I_i \right)\;,\\
\end{aligned}
\end{equation*} 
is the cumulative distribution function of the couple 
$(m_{n+1}(\x_{n+1}^q(B_i)-G(x)), Z):=(W, Z)\;,$ at point $(0, I_i)$. 
This random vector, conditionally on $\mathcal{A}_n$ is Gaussian. We denote by $M$ and $R$ its mean vector and covariance matrix, respectively. 

Thanks to Proposition \ref{MAJkrigeage}, we have:
\begin{equation}
\label{MAJpointquantile}
\begin{aligned}
m_{n+1}(\x_{n+1}^q(B_i))=m_n(\x_{n+1}^q(B_i))-k_n(\x_{n+1}^q(B_i), \x_{n+1})Z\;,\\
\end{aligned}
\end{equation}
which gives
$$M=
\begin{pmatrix}
m_n(\x_{n+1}^q(B_i)) - m_n(x) \\
0\\
\end{pmatrix}, \, \, R= \begin{pmatrix} \Var(W) & \Cov(W, Z)\\
\Cov(W, Z) & \Var(Z) \\
\end{pmatrix}\;,$$
with
$$\Var(W):=\sigma_{W}=s_n(x)^2 + \frac{k_n(\x_{n+1}^q(B_i), \x_{n+1})^2}{s_n(\x_{n+1})^2} - 2 \frac{k_n(\x_{n+1}^q(B_i), \x_{n+1})k_n(x, \x_{n+1})}{s_n(\x_{n+1})^2},$$

$$\Cov(W, Z)=\frac{k_n(\x_{n+1}^q(B_i), \x_{n+1})-k_n(x, \x_{n+1})}{s_n(\x_{n+1})^2}\, \, \text{and} \, \, \Var(Z)=\frac{1}{s_n(\x_{n+1})^2}\;.$$

We can conclude by centering and normalizing:
\begin{equation*}
\begin{aligned}
T_n&  = \P_n \left( W  \leq 0 \cap Z \leq I_i \right)\\
&= \P \left( \frac{W-(m_n(\x_{n+1}^q(B_i)) - m_n(x) )}{\sqrt{Var(W)}} \leq \frac{m_n(x)-m_n(\x_{n+1}^q(B_i))}{\sqrt{Var(W)}} \cap s_n(\x_{n+1})Z \leq I_i s_n(\x_{n+1}) \right)\\
& := \P \left(S \leq e^i_n(\x_{n+1}; x; \x_{n+1}^q(B_i)) \cap T \leq f^i_n(\x_{n+1};I_i) \right)\;,\\
\end{aligned}
\end{equation*} 

where $(S, T)$ is a Gaussian random vector of law $\mathcal{N} \left( 0, \begin{pmatrix} 1 & r^i_n \\
r^i_n & 1 \\
\end{pmatrix}\; \right)$ with
$$r^i_n:=r_n(\x_{n+1}; x; \x_{n+1}^q(B_i))=\frac{k_n(\x_{n+1}^q(B_i), \x_{n+1})- k_n(x, \x_{n+1})}{\sqrt{\Var(W)} s_n(\x_{n+1})}\;.$$ 
%
%
$$e^i_n(\x_{n+1}; x; \x_{n+1}^q(B_i))= \frac{m_n(x) - m_n(\x_{n+1}^q(B_i))}{\sqrt{\Var(W)}}\;,$$ and $$f^i_n(\x_{n+1}; I_i)= I_is_n(\x_{n+1}).$$

Finally, we get for $1 \leq i \leq L$,
$$\P_n(G(x) \geq m_{n+1}(\x_{n+1}^q(B_i)) \cap Z \leq I_i)= \Phi_{r^i_n}\left(e^i_n(\x_{n+1}; x; \x_{n+1}^q(B_i)), f^i_n(\x_{n+1}; I_i)\right)\;,$$
where we denote $\Phi_r$ the cumulative distribution function of the centered Gaussian random vector of covariance matrix 

$$\Sigma=\begin{pmatrix} 1 & r \\
r & 1 \\
\end{pmatrix}\;.$$

Similarly, for $0 \leq i \leq L$:
$$\P_n(G(x) \geq m_{n+1}(\x_{n+1}^q(B_i)) \cap Z \leq I_{i+1})= \Phi_{r^i_n}\left(e^i_n(\x_{n+1}; x; \x_{n+1}^q(B_i)), f^i_n(\x_{n+1}; I_{i+1})\right)\;,$$
and 
$$\P_n(G(x) \geq m_{n+1}(\x_{n+1}^q(B_i)) \cap Z \geq I_{L})= \Phi_{-r^i_n}\left(e^i_n(\x_{n+1}; x; \x_{n+1}^q(B_i)), -f^i_n(\x_{n+1}; I_{L})\right)\;.$$

%
%
%
\noindent {\bf 6.2. Proof of Proposition \ref{var}}\label{sec:proofvar}

We first recall the following total variance law formula:
\begin{lemma}
Let $E_1$, $\dots$ $E_n$ be mutually exclusive and exhaustive events. Then, for a random variable $U$, the following equality holds:
\begin{equation*}
\label{variancelaw}
\begin{aligned}
Var(U) &= \sum_{i=1}^{n}{\Var(U \mid E_i) \P(E_i)}+ \sum_{i=1}^n {\E(U \mid E_i)^2 (1-\P(E_i))\P(E_i)}\\
	&- 2\sum_{i=2}^{n}\sum_{j=1}^{i-1}\E(U \mid E_i) \P(E_i)\E(U \mid E_j) \P(E_j)\;.
\end{aligned}
\end{equation*} 
\end{lemma}

In our case, we want to compute $\Var(q_{n+1}(\x_{n+1})|\mathcal{A}_n):=\Var_n(q_{n+1}(\x_{n+1}))$. 
Since the events $\{Z\in  B_i\}_{ 1 \leq i \leq L}$ are mutually exclusive and exhaustive, we can apply Lemma \ref{variancelaw}:
\begin{equation*}
\begin{aligned}
\Var_n(q_{n+1}(\x_{n+1}))&= \displaystyle\sum_{i=1}^L \Var_n(m_{n+1}(\x_{n+1}^q(B_i)) | Z \in B_i) \P_n(Z \in B_i)\\
& + \displaystyle\sum_{i=1}^L \E_n\left(m_{n+1}(\x_{n+1}^q(B_i))| Z \in B_i\right)^2 \left( 1 - \P_n(Z \in B_i) \right) \P_n(Z \in B_i)\\
& - 2 \displaystyle\sum_{i=2}^L \displaystyle\sum_{j=1}^{i-1} \E_n\left( m_{n+1}(\x_{n+1}^q(B_i)) | Z \in B_i \right) \P_n \left( Z \in B_i\right)\\
& \times \E_n\left( m_{n+1}(\x_{n+1}^q(B_j)) | Z \in B_j \right) \P_n \left( Z \in B_j\right)\;.\\
\end{aligned}
\end{equation*}

Thanks to equation (\ref{MAJpointquantile}), we get
\begin{equation*}
\begin{aligned}
m_{n+1}(\x_{n+1}^q(B_i))=m_n(\x_{n+1}^q(B_i))-k_n(\x_{n+1}^q(B_i), \x_{n+1})Z\;.\\
\end{aligned}
\end{equation*}

Then,
\begin{equation*}
\begin{aligned}
\Var_n(q_{n+1}(\x_{n+1}))&= \displaystyle\sum_{i=1}^n k_n(\x_{n+1}^q(B_i), \x_{n+1})^2 \Var_n(Z| Z \in B_i) \P_n(Z \in B_i)\\
& + \displaystyle\sum_{i=1}^L \left( m_n(\x_{n+1}^q(B_i)) - k_n(\x_{n+1}^q(B_i), \x_{n+1})\E_n\left(Z| Z \in B_i\right)\right)^2 \\
& \times \left( 1 - \P_n(Z \in B_i) \right) \P_n(Z \in B_i)\\
& - 2 \displaystyle\sum_{i=2}^L \displaystyle\sum_{j=1}^{i-1} \left(m_n(\x_{n+1}^q(B_i)) - k_n(\x_{n+1}^q(B_j), \x_{n+1})\E_n\left(Z | Z \in B_j \right)\right) \P_n \left( Z \in B_i\right)\\
& \times \left(m_n(\x_{n+1}^q(B_j)) - k_n(\x_{n+1}^q(B_j), \x_{n+1})\E_n\left(Z | Z \in B_j \right)\right) \P_n \left( Z \in B_j\right)\;.\\
\end{aligned}
\end{equation*}

Since $Z$ is a centered Gaussian random variable of variance $s_n(\x_{n+1})^{-2}$  we have:
$$P_i:= \P_n(Z \in B_i)= \Phi(s_n(\x_{n+1})I_{i+1})-\Phi(s_n(\x_{n+1})I_i)\;.$$

To conclude, we have now to find analytical forms for the quantities $\Var \left(Z \mid I_i < Z < I_{i+1} \right)$ and $\E \left(Z \mid I_i < Z < I_{i+1} \right)$. 
To do so, let us use the following result on truncated Gaussian random variable (see \cite{talis} for proofs): 
\begin{lemma}
\label{truncated}
Let $U$ be a real random variable such that $U \sim \mathcal{N}(\mu, \sigma)$. Let $u$ and $v$ be two real numbers. 
We have:
\begin{equation*}
\begin{aligned}
\E(U | u <U<v)&= \mu + \frac{\phi \left( \frac{v- \mu}{\sigma} \right) - \phi \left( \frac{u- \mu}{\sigma} \right)}{\Phi \left( \frac{v-\mu}{\sigma} \right) - \Phi \left( \frac{w - \mu}{\sigma} \right)} \sigma \;,\\
\end{aligned}
\end{equation*}

\begin{equation*}
\begin{aligned}
\Var(U | u<U<v)= \sigma^2 \left[1 + 
\frac{\frac{u-\mu}{\sigma} \phi \left(\frac{u-\mu}{\sigma}\right) - \frac{v-\mu}{\sigma} \phi \left(\frac{v-\mu}{\sigma}\right)}{\Phi \left(\frac{v-\mu}{\sigma} \right) - \Phi \left( \frac{u-\mu}{\sigma} \right)} 
- \left( \frac{\phi \left(\frac{u-\mu}{\sigma}\right)- \phi \left(\frac{v-\mu}{\sigma}\right)}{\Phi \left(\frac{v-\mu}{\sigma} \right) - \Phi \left( 
\frac{u-\mu}{\sigma} \right)} \right)^2 \right]\;.\\
\end{aligned}
\end{equation*}
\end{lemma}

We apply Lemma \ref{truncated} for $U=Z$, $u=I_i$ and $v=I_{i+1}$ and conclude that

$$E(s_n(\x_{n+1}), I_{i+1}, I_i):=E_n(Z|Z \in B_i)= \frac{\phi(s_n(\x_{n+1})I_{i})-\phi(s_n(\x_{n+1})I_{i+1})}{\Phi(s_n(\x_{n+1})I_{i+1})- \Phi(s_n(\x_{n+1})I_i)} \frac{1}{s_n(\x_{n+1})}\;,$$

\begin{equation*}
\begin{aligned}
V(s_n(\x_{n+1}), I_{i+1}, I_i):&=\Var_n(Z | Z \in B_i)\\
&= \frac{1}{s_n(\x_{n+1})^2} \left[ 1+ \frac{I_is_n(\x_{n+1}) \phi(s_n(\x_{n+1}) I_{i}) - s_n(\x_{n+1})I_{i+1} \phi(s_n(\x_{n+1}) I_{i})}{\Phi(s_n(\x_{n+1}) I_{i+1})-\Phi( s_n(\x_{n+1})I_{i})} \right. \\
&\left. - \left( \frac{\phi(s_n(\x_{n+1})I_{i}) - \phi(s_n(\x_{n+1})I_{i+1})}{\Phi(s_n(\x_{n+1})I_{i+1}) - \Phi(s_n(\x_{n+1})I_i)} \right)^2\right]\;.\\
\end{aligned}
\end{equation*}

%
%

\noindent {\bf 6.3 Test functions}\label{sec:test}

Two-dimensional Branin function:
\begin{equation}\label{eq:branin}
  g(\mathbf{x}) = \left( \bar{x}_2 - \frac{5.1 \bar{x}_1^2}{4 \pi^2} + \frac{5 \bar{x}_1}{\pi} - 6 \right)^2
+ \left( 10 - \frac{10}{8 \pi} \right) \cos(\bar{{x}}_1) +10
\end{equation}
with: $\bar{x}_1 = 15 \times x_1-5$, $\bar{x}_2 = 15 \times x_2$.

Four-dimensional Hartman function:
\begin{eqnarray}\label{eq:hartman}
 g(\x) = \frac{-1}{1.94} \left[2.58 + \sum_{i=1}^4{ C_i \exp \left( -\sum_{j=1}^4{  a_{ji} \left( x_j - p_{ji} \right)^2  } \right) } \right],
\end{eqnarray}
with
\begin{small}
\begin{eqnarray*}
\mathbf{C} = \left[ \begin{array}{c} 
1.0\\ 1.2\\ 3.0\\ 3.2 
\end{array} \right], & 
\mathbf{a} = \left[ \begin{array}{cccc}
	         10.00 &  0.05&  3.00& 17.00 \\
           3.00& 10.00&  3.50&  8.00\\
           17.00& 17.00&  1.70&  0.05\\
           3.50&  0.10& 10.00& 10.00
\end{array} \right],
\mathbf{p} = \left[ \begin{array}{cccc}
           0.1312& 0.2329& 0.2348& 0.4047\\
           0.1696& 0.4135& 0.1451& 0.8828\\
           0.5569& 0.8307& 0.3522& 0.8732\\
           0.0124& 0.3736& 0.2883& 0.5743
\end{array} \right].
\end{eqnarray*}
\end{small}

Six-dimensional Ackley function:
\begin{equation}\label{eq:ackley}
 g(\x) = 20 + \exp(1) -20 \exp \left( -0.2 \sqrt{\frac{1}{4}\sum_{i=1}^4 x_i^2} \right)  - \exp \left[\frac{1}{4}\sum_{i=1}^4 \cos \left(2\pi x_i \right) \right]
\end{equation}

\bibliographystyle{plain}
\bibliography{bib}

\begin{thebibliography}{10}

\bibitem{flooding}
Aur{\'e}lie Arnaud, Julien Bect, Mathieu Couplet, Alberto Pasanisi, and
  Emmanuel Vazquez.
\newblock {{\'E}valuation d'un risque d'inondation fluviale par planification
  s{\'e}quentielle d'exp{\'e}riences}.
\newblock In {\em {42{\`e}mes Journ{\'e}es de Statistique}}, Marseille, France,
  France, 2010.

\bibitem{bect}
Julien Bect, David Ginsbourger, Ling Li, Victor Picheny, and Emmanuel Vazquez.
\newblock Sequential design of computer experiments for the estimation of a
  probability of failure.
\newblock {\em Statistics and Computing}, 22(3):773--793, 2012.

\bibitem{cannamela}
Claire Cannamela, Josselin Garnier, and Bertrand Iooss.
\newblock Controlled stratification for quantile estimation.
\newblock {\em The Annals of Applied Statistics}, pages 1554--1580, 2008.

\bibitem{chevalierExcursion}
Cl{\'e}ment Chevalier, David Ginsbourger, Julien Bect, Emmanuel Vazquez, Victor
  Picheny, and Yann Richet.
\newblock Fast parallel kriging-based stepwise uncertainty reduction with
  application to the identification of an excursion set.
\newblock {\em Technometrics}, 56(4), 2014.

\bibitem{MAJkrigeage}
Cl{\'e}ment Chevalier, David Ginsbourger, and Xavier Emery.
\newblock Corrected kriging update formulae for batch-sequential data
  assimilation.
\newblock In {\em Mathematics of Planet Earth}, pages 119--122. Springer, 2014.

\bibitem{order}
David and Nagaraja.
\newblock {\em Order Statistics}.
\newblock Wiley, 2003.

\bibitem{dixon1978towards}
L.C.W. Dixon and G.P. Szeg{\"o}.
\newblock {\em Towards Global Optimisation 2}, volume~2.
\newblock North Holland, 1978.

\bibitem{regret}
Nando~D. Freitas, Masrour Zoghi, and Alex~J. Smola.
\newblock Exponential regret bounds for gaussian process bandits with
  deterministic observations.
\newblock In John Langford and Joelle Pineau, editors, {\em Proceedings of the
  29th International Conference on Machine Learning (ICML-12)}, pages
  1743--1750, New York, NY, USA, 2012. ACM.

\bibitem{geman1996active}
Donald Geman and Bruno Jedynak.
\newblock An active testing model for tracking roads in satellite images.
\newblock {\em Pattern Analysis and Machine Intelligence, IEEE Transactions
  on}, 18(1):1--14, 1996.

\bibitem{gramacy2008bayesian}
Robert~B Gramacy and Herbert~KH Lee.
\newblock Bayesian treed gaussian process models with an application to
  computer modeling.
\newblock {\em Journal of the American Statistical Association}, 103(483),
  2008.

\bibitem{jala1}
Marjorie Jala, C{\'e}line L{\'e}vy-Leduc, Eric Moulines, Emmanuelle Conil, and
  Joe Wiart.
\newblock Sequential design of computer experiments for parameter estimation
  with application to numerical dosimetry.
\newblock In {\em Signal Processing Conference (EUSIPCO), 2012 Proceedings of
  the 20th European}, pages 909--913. IEEE, 2012.

\bibitem{jala2}
Marjorie Jala, C{\'e}line L{\'e}vy-Leduc, {\'E}ric Moulines, Emmanuelle Conil,
  and Joe Wiart.
\newblock Sequential design of computer experiments for the assessment of fetus
  exposure to electromagnetic fields.
\newblock {\em Technometrics}, 58(1):30--42, 2016.

\bibitem{kenkel2011pbivnorm}
B~Kenkel.
\newblock {\em pbivnorm: Vectorized Bivariate Normal CDF}, 2012.
\newblock R package version 0.5-1.

\bibitem{kennedy2001bayesian}
Marc~C Kennedy and Anthony O'Hagan.
\newblock Bayesian calibration of computer models.
\newblock {\em Journal of the Royal Statistical Society: Series B (Statistical
  Methodology)}, 63(3):425--464, 2001.

\bibitem{algo}
Dong~C Liu and Jorge Nocedal.
\newblock On the limited memory bfgs method for large scale optimization.
\newblock {\em Mathematical programming}, 45(1-3):503--528, 1989.

\bibitem{oakley}
Jeremy Oakley.
\newblock Estimating percentiles of uncertain computer code outputs.
\newblock {\em Journal of the Royal Statistical Society: Series C (Applied
  Statistics)}, 53(1):83--93, 2004.

\bibitem{victorMO}
Victor Picheny.
\newblock Multiobjective optimization using gaussian process emulators via
  stepwise uncertainty reduction.
\newblock {\em Statistics and Computing}, pages 1--16, 2013.

\bibitem{rasmussen2006gaussian}
C.E. Rasmussen and C.K.I. Williams.
\newblock {\em {Gaussian processes for machine learning}}.
\newblock MIT Press, 2006.

\bibitem{roustant2012dicekriging}
Olivier Roustant, David Ginsbourger, and Yves Deville.
\newblock Dicekriging, diceoptim: Two r packages for the analysis of computer
  experiments by kriging-based metamodeling and optimization.
\newblock {\em Journal of Statistical Software}, 51(1):1--55, 2012.

\bibitem{sacks1989design}
Jerome Sacks, William~J Welch, Toby~J Mitchell, and Henry~P Wynn.
\newblock Design and analysis of computer experiments.
\newblock {\em Statistical science}, pages 409--423, 1989.

\bibitem{covfunction}
Michael~L. Stein.
\newblock {\em Interpolation of spatial data: some theory for kriging}.
\newblock Springer Science \& Business Media, 2012.

\bibitem{talis}
Georges~M. Tallis.
\newblock The moment generating function of the truncated multi-normal
  distribution.
\newblock {\em Journal of the Royal Statistical Society. Series B
  (Methodological)}, pages 223--229, 1961.

\bibitem{villemonteix2009informational}
Julien Villemonteix, Emmanuel Vazquez, and Eric Walter.
\newblock An informational approach to the global optimization of
  expensive-to-evaluate functions.
\newblock {\em Journal of Global Optimization}, 44(4):509--534, 2009.

\end{thebibliography}

\end{document}